\documentclass[12pt,letterpaper,usenames,dvipsnames]{article}
\usepackage{jheppub}
\makeatletter
\def\@fpheader{\relax \phantom{.}}
\makeatother

\allowdisplaybreaks[2]        
\usepackage{bm}                
\usepackage{bbold}            
\usepackage{amscd}           

\usepackage{tikz}                  
\usetikzlibrary{matrix}
\usetikzlibrary{arrows, shapes.arrows}
\usetikzlibrary{decorations, decorations.markings}

\usepackage{xcolor}   

\usepackage{url}

\def\nn{\nonumber}
\def\bar{\overline}
\def\til{\widetilde}

\def\vev#1{{\langle{#1}\rangle}} 
\def\I{\mathbb{1}}

\newcommand{\minus}{\setminus}

\def\kn{K^n_{p,q}}
\def\kone{K_{p,q}}
\newcommand{\link}[1]{K^{n_{#1}}_{p_{#1},q_{#1}}}
\def\C{\mathbb{C}} 
\def\N{\mathbb{N}}

\def\R{\mathbb{R}} 
\def\Z{\mathbb{Z}} 
\def\bl{{\boldsymbol{\ell}}}
\def\fL{{\mathfrak L}}
\def\bm{{\bf m}}
\def\cN{{\mathcal N}}
\def\a{{\alpha}}
\def\ta{{\til\a}}
\def\b{{\beta}}
\def\tb{{\til\b}}
\def\g{{\gamma}}

\def\d{{\delta}}
\def\th{{\theta}}
\def\ps{{\psi}}

\title{Fundamental groups for torus link complements}
\author[a]{Philip C. Argyres,}
\author[b]{Dnyanesh P. Kulkarni}
\affiliation[a]{Department of Physics, University of Cincinnati,
Cincinnati, OH, USA 45221}
\affiliation[b]{Department of Physics, Cornell University, Ithaca, NY, USA 14853}
\emailAdd{philip.argyres@gmail.com}
\emailAdd{dpk52@cornell.edu}

\abstract{For an arbitrary positive integer $n$ and a pair $(p, q)$ of coprime integers, consider $n$ copies of a torus $(p,q)$ knot placed parallel to each other on the surface of the corresponding auxiliary torus: we call this assembly a torus $n$-link.  We compute economical presentations of knot groups for torus links using the groupoid version of the Seifert--van Kampen theorem.  Moreover, the result for an individual torus $n$-link is generalized to the case of multiple ``nested" torus links, where we inductively include a torus link in the interior (or the exterior) of the auxiliary torus corresponding to the previous link.  The results presented here have been useful in the physics context of classifying moduli space geometries of four-dimensional ${\mathcal N}=2$ superconformal field theories.}
 
\begin{document}
\maketitle

\section{Introduction}
\label{sec:intro}

\begin{figure}[ht]
\centering
\includegraphics[width=.50 \textwidth]{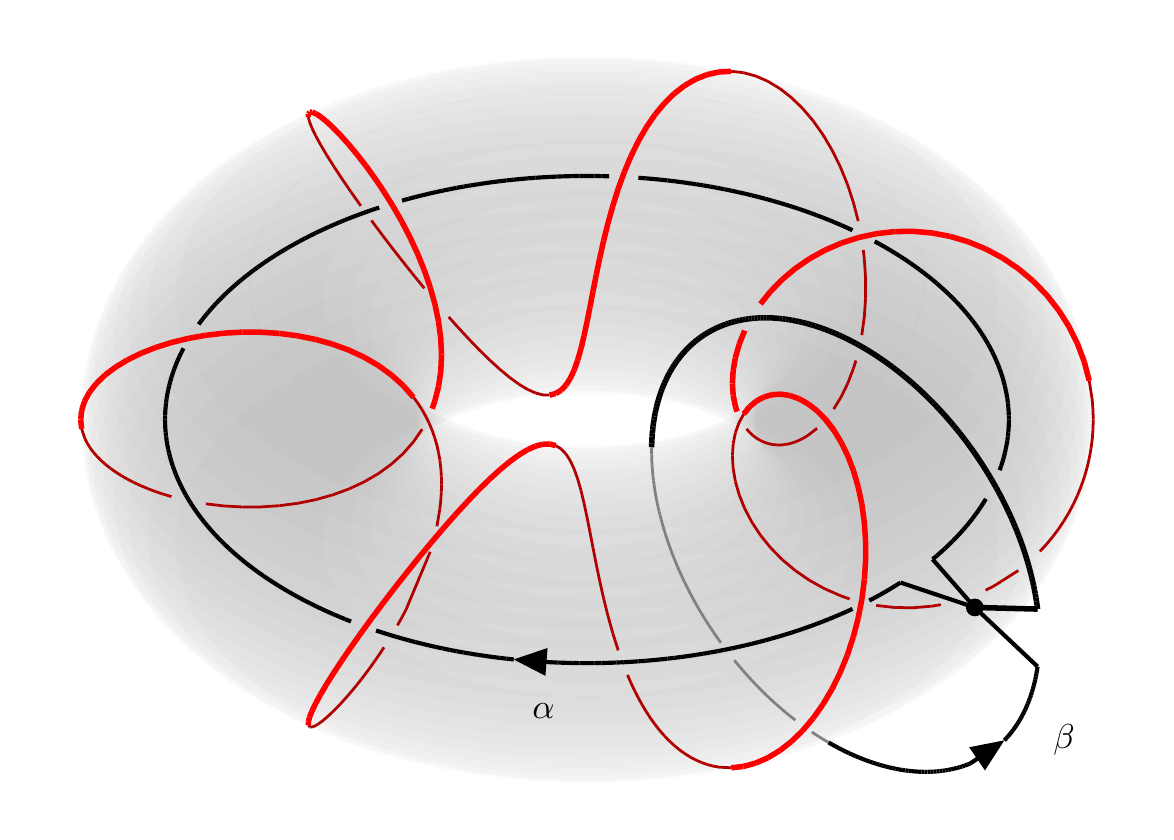}
\caption{Depiction of a $K_{1,6}$ torus knot (in this case homeomorphic to an unknot in $S^3$) consisting of the red circle.  The gray donut (solid torus) is there for visualization purposes; it is the set $D_A$ introduced in section \ref{sec:definitions}, and its boundary is the auxiliary 2-torus $T$.  The $\a$ cycle threads the interior of the donut, while $\b$ threads the hole of the donut.}
\label{fig1}
\end{figure}

A knot is a smooth embedding of $S^1$ into $\R^3$ or $S^3$ up to equivalence by ambient isotopies.  Let $\kone$ denote a torus knot for $p$, $q$ coprime integers.  $\kone$ can be visualized as a simple closed curve lying on the surface of an auxiliary 2-torus, $T$, embedded in $\R^3$.  The meridian cycle, $\bm\in H_1(T)$, of the torus is a cycle which is a primitive homology generator which can be deformation retracted to a point in the interior of $T$, and assign $T$ 
an orientation.  The longitudinal cycle, $\bl\in H_1(T)$, is the cycle with intersection $\vev{\bl,\bm}=1$ relative to some chosen orientation on $T$.  Then $\kone$ is a simple curve that winds $q$ times around the longitude while it winds $p$ times around the meridian of $T$, and thus its homology class is $[\kone]=p\bl+q\bm$.  See figure \ref{fig1}.  

There are restrictions on the pairs of integers $(p,q)$ describing torus knots.  One is a convention that $(p,q)\neq(0,0)$, since otherwise it is retractable to a point within $T$ so not considered to be a torus knot.  Another restriction, mentioned above,  is that $\gcd(p,q)=1$, for otherwise a simple closed curve in $T$ with these winding numbers does not exist:  to avoid self crossings it would have to be a link with $n=\gcd(p,q) >1$ components.  

If one is only interested in knots up to orientation-preserving homeomorphisms, then there would also be the identifications $\kone \simeq K_{q,-p}$ and $K_{p,1} \simeq K_{0,1} \simeq K_{0,-1}$.  ($K_{0,1}$ is an ``unknot''.)  Although our goal is to compute knot groups --- which are topological invariants --- it will be convenient not to impose the above identifications in our notation describing torus knots and related classes of links which we will duly introduce.  

The knot group of $\kone$ is the fundamental group of the path-connected space $\R^3 \minus \kone$.  This knot group can be computed using the Seifert--van Kampen theorem, and a presentation for it in terms of generators and relations is 
\begin{align}\label{centuryold}
\pi_1(\R^3 \minus \kone) = \vev{\a, \b \, | \, \a^p \b^{-q} } .
\end{align}
See, e.g., example 1.24 in \cite{Hatcher}.  Given a choice of base point, cycles corresponding to the generators $\a$ and $\b$ are shown in figure \ref{fig1}.  In the case of an unknot, $(p,q) = (1,0)$, the fundamental group is $\Z$.  It is also easy to check from this presentation that the fundamental group is $\Z$ also for $(p,q)=(1,q)$ for any $q$, as must be the case since the knots are homeomorphic.

The knot group of $\kone$ can also be computed in the Wirtinger presentation \cite{Tietze} by suitably projecting the knot to a plane, keeping track of over- or under-crossings, assigning generators as cycles looping around each segment of the projected knot between crossings, and imposing suitable relations at each crossing.  Though equivalent to the Seifert--van Kampen presentation, the Wirtinger presentation is much less economical, as even after ``immediate" simplifications it involves $p$ generators and $q-1$ relations.

The main construction of this paper generalizes the above situation to \emph{torus links}.  A link is an embedding of a finite disjoint union of circles into $S^3$.  For a positive integer $n$, consider $n$ copies of $\kone$ situated parallel to each other on $T$.  See figure \ref{fig3} for an illustration.  We denote the union of these $n$ knots as $\kn$, and call this assembly a torus $(p,q)$ $n$-link.  (So $\kone = K^1_{p,q}$.)  Our goal is to compute the knot group of torus links in an economical presentation like that of the Seifert--van Kampen presentation for torus knots.

The problem is that a direct application the Seifert--van Kampen theorem fails for torus links because the intersection of two spaces used to cover $\R^3 \minus \kn$ is not path connected for $n>1$.  However, there exists a beautiful way of conquering such difficulties using the notion of fundamental groupoids.  The basic idea is to select multiple base points --- one in each of the $n$ connected components of the intersection, which is $T \minus \kn$ --- and then construct fundamental groupoids with respect to these base points for each of the spaces in a Seifert--van Kampen type construction, to obtain the fundamental groupoid for the space of interest, $\R^3 \minus \kn$.  One then reduces (uniquely) the fundamental groupoid to the fundamental group as a last step.  The computation is described in section \ref{sec:links}, where the result is given in \nameref{result1}.

Knowing the fundamental groupoid allows us to generalize the problem further: given a torus link $\link{1}$, we can include a torus link $\link{2}$ in the interior (or the exterior) of the corresponding auxiliary torus $T_1$, and another link $\link{3}$ in the interior (exterior) of the second auxiliary torus $T_2$, and so on.  We establish here an inductive argument to carry out the explicit computation of the fundamental groups of the complements in $\R^3$ of such nested torus links. The result is given in \nameref{nstdgrp}.

The rest of this paper makes all the above remarks more precise, introduces the tools to be used, and applies them to explicit computations of the fundamental groups.

We finally remark that the original motivation for this problem comes from the physics of certain superconformal quantum field theories.  Specifically, understanding the topology of nested torus links turns out to be useful for constructing the special K\"ahler \cite{Freed} Coulomb branch geometries (a.k.a., Seiberg-Witten curves \cite{SW}) of 4-dimensional $\cN=2$ superconformal quantum field theories when the rank of the gauge group is 2.  An application of the results found here to this physical problem has appeared in \cite{Argyres}.

\section{Preliminaries and definitions}
\label{sec:definitions}

We can consider torus knots or links as embedded in either $\R^3$ or its one-point compactification $S^3$.  The fundamental group of the complement of the knot or link is the same under either choice of the ambient space; see, e.g., example 1.24 in \cite{Hatcher}.  It will be convenient to take $S^3$ to be the ambient space.  Although in $S^3$ there is no longer an invariant distinction between the interior and the exterior of the torus $T$, we will nevertheless conventionally call one connected component, $D_A$, of $S^3 \minus T$ the ``interior", and the other, $D_B$, the ``exterior".   (In general, we will use the letter $A$ to denote sets interior to $T$ and $B$ those exterior to $T$ in this convention.)

We begin by giving explicit parametrizations of the various sets, $S^3$, $D_A$, $D_B$, $T$, $\kn$, $A$, $B$, $S^1_A$, and $S^1_B$, that we will use.  These parameterization are useful for making explicit the  deformation retractions we use below.  Describe $S^3$ as the set of solutions to $|u|^2 + |v|^2 = 2$ with $(u,v)\in\C^2$.  These can be parameterized by $(u,v)= \sqrt2 (\cos\chi \, e^{i\th}, \sin\chi\, e^{i\ps})$, with $0\le \chi \le \pi/2$, and $\th$ and $\ps$ both periodic with period $2\pi$.  $D_A = \{ (\chi,\th,\ps)\, |\, \chi<\pi/4\}$, $D_B = \{ (\chi,\th,\ps)\, |\, \chi>\pi/4\}$, and the 2-torus which is the boundary of both of these sets is $T= \{ (\chi,\th,\ps)\, |\, \chi=\pi/4\}$.   Thus $D_A$ and $D_B$ are donuts (solid tori).  These can also be described as the intersection of sets in $\C^2$ with $S^3$ as:  $D_A =\{ |v|<1\} \cap S^3$, $T =\{ |u|=|v|=1\} \cap S^3$ and $D_B =\{ |u|<1\} \cap S^3$.  The basic homotopies we will need in the arguments to follow simply involve shrinking or expanding the tori defined by constant $\chi$, thus involve only shifts and rescalings of the $\chi$ coordinate. 

The $\kn$ torus link is $n$ parallel copies of a knot placed on the $\chi=\frac\pi4$ torus which winds $q$ times around the $\th$ circle as is winds $p$ times around the $\ps$ circle,
\begin{align}\label{Knqpdef}
\kn = \left\{ (\chi,\th,\ps)\, \Big|\, \chi=\frac\pi4 \, ,\   
q\th - p\ps = 2\pi \frac kn \mod 2\pi \, , \ \text{and}\ k\in\Z_n 
\right\} ,
\end{align}  
where $k$ labels the $n$ components of the link.  Alternatively, $\kn = \{ u^{nq}=v^{np}\}\cap T$.  

It will also be useful to define the circles $S^1_A = \{ (\chi,\th) | \chi=0\}$ and $S^1_B = \{ (\chi,\ps) | \chi=\pi/2\}$ (the coordinate system degenerates at these two circles, so $\ps$ is redundant at $S^1_A$ and $\th$ at $S^1_B$).  Alternatively, these can be described as the intersections $S^1_A=\{ |v|=0\} \cap S^3$ and $S^1_B=\{ |u|=0\} \cap S^3$.   $D_A$ is homotopic to the interior circle $S^1_A$ by shrinking the $0\le\chi\le\pi/4$ interval to $\chi=0$.  Likewise, $D_B$ is homotopic to the exterior circle $S^1_B$ shrinking the  $\pi/4\le\chi\le\pi/2$ interval to $\chi=\pi/2$.

All these sets are also easily visualized:  the figures give examples where the $S^3$ is stereographically projected to $\R^3$ with $(\chi,\th)=(0,0)$ is projected to the origin of $\R^3$, and $(\chi,\psi)=(\pi/2,0)$ is projected to infinity.   Then the interior of $T$ is the shaded solid donut of the figures.  For example, figure \ref{fig1} shows a $K_{1,6}$ torus knot, figure \ref{fig3} shows a $K^3_{1,3}$ torus link, and figure \ref{fig4} shows a $K_{1,6}$ torus knot linked with two unknots placed at $S^1_A$ and $S^1_B$.

Next, we review how to use the Seifert--van Kampen theorem \cite{Seifert, vanKampen} to compute the knot group of a torus knot.  Consider the torus knot $\kone$ that resides on the auxiliary torus $T$.  Then define $A := D_A \cup (T \minus \kone)$ and $B := D_B \cup (T \minus \kone)$, so that $C := A \cap B = T \minus \kone$ and $A \cup B = S^3 \minus \kone$.  Each of these spaces is path-connected.  In particular, $A$, $B$, and $C$ are each deformation retractable to an $S^1$ so $\pi_1(A) =\vev{\a \, | \, \varnothing}$, $\pi_1(B) = \vev{\b \, | \, \varnothing}$, and $\pi_1(C) = \vev{\g \, | \, \varnothing}$.  Here we are describing the fundamental groups as $\vev{\text{generators}\,|\,\text{relations}}$; we use $\varnothing$, the empty set, to denote no relations.   Since $C$ is a subset of both $A$ and $B$, there are homomorphisms $i:\pi_1(C)\to\pi_1(A)$ and $j:\pi_1(C)\to\pi_1(A)$ from \lq\lq{}pushing out\rq\rq{} the homotopy classes of closed paths in $C$ to those in $A$ or $B$ using the $C\subset A$ and $C\subset B$ inclusion maps.  Furthermore, it is easy to see from our explicit coordinatization of the torus that we can choose $i(\g) = \a^p$ and $j(\g)=\b^q$.  The Seifert--van Kampen theorem then implies that $ \pi_1 (S^3 \minus \kone) = \pi_1 (A) \ast_{\pi_1(C)} \pi_1(B) = \vev{\a, \b \, | \, \a^ p = \b^ q}$, which is \eqref{centuryold}.  The $\ast_{\pi(C)}$ notation denotes the pushout, or amalgamated free product, of groups.  It is the free product with additional relations $i(\g)=j(\g)$ for all generators $\g\in\pi_1(C)$.

Now, consider a similar construction for a torus link $\kn$. That is, we have 
\begin{align}\label{Knpqsets}
A &:= D_A \cup (T \minus \kn), & 
B &:= D_B \cup (T \minus \kn), \nn\\
C &:= A \cap B = T \minus \kn, & 
A \cup B &= S^3 \minus \kn.
\end{align}
But the space $C$ is not path-connected, and the Seifert--van Kampen theorem cannot be applied to compute the knot group of $\kn$.  

Note, however, that $C$ has $n$ distinct path components each of which deformation retracts onto a homeomorphic image of knot $\kone$. 
Informally, each of these path components deformation retracts to the line that runs through its center, and since this middle line runs parallel to the $n$ connected components of the link $\kn$, it is homeomorphic to $\kone$.  Since the spaces $A, B, A\cup B$ are all path-connected and $C$ is symmetric with respect to its $n$ distinct path components (i.e., the distinct path components are pairwise homeomorphic), one suspects that there is some simple way of generalizing the Seifert--van Kampen argument.   Such a way is provided by the notion of a fundamental groupoid and by the groupoid generalization of the Seifert--van Kampen theorem \cite{Brown67}.  

The rest of this section gives a quick review of the definition and properties of fundamental groupoids and their presentations in terms of generating graphs and relations.  More detailed expositions can be found in \cite{May, Brown06}.

Groupoids are small categories in which every morphism is an isomorphism.  It is a standard terminology to refer to the isomorphisms in a groupoid as its \emph{elements}.  The fundamental groupoid $\pi_1(X,P)$ of a space $X$ with respect to a set of base points $P\subset X$ is the groupoid whose objects are points in $P$, and whose elements are homotopy equivalence classes of paths from $x$ to $y$ with $x, y \in P$.

We will only use a \emph{finite pointed} model for fundamental groupoids, where we take the object set $P$ to be finite.  If set $P$ has only one element, $P = \{x\}$, the fundamental groupoid is the same as the fundamental group defined with the given base point, i.e.,  $\pi_1(X,P) = \pi_1(X,x)$.  In general, the elements of $\pi_1 (X, P)$ --- the homotopy equivalence classes in $X$ of paths from $x$ to $y$ in $P$ --- induce respective isomorphisms between the fundamental groups $\pi_1(X,x)$ and $\pi_1(X,y)$ for different base points.

The object set $P$ is \emph{representative} if $P$ non-trivially intersects each path component of $X$.  All the groupoids that we use are defined over representative subsets of the corresponding topological spaces.

A groupoid is \emph{connected} if there is an isomorphism between every pair of its objects.  It then follows that the fundamental groupoid defined over a representative set of a topological space is connected if and only if the space itself is path-connected.  Since a groupoid with exactly one object is a group, the full subgroupoid of a groupoid defined at a chosen object is the corresponding \emph{object group}.  An object group of a fundamental groupoid is the fundamental group of the space with respect to the chosen object, i.e., the base point.

A representative base point set $P \subset X$ will be chosen and remain fixed for each space we discuss, so we will use a simplified notation,
\begin{align}\label{Notn}
\pi X : = \pi_1 (X, P),
\end{align}
for the corresponding fundamental groupoid of a space $X$.   The set of elements (isomorphisms) of $X$ connecting $x$ to $y$ for $x,y\in P$ is denoted $\pi X(x,y)$.  The object group at $x$ is denoted $\pi X(x) := \pi X(x,x)$, so 
\begin{align}\label{Notn2}
\pi X(x) = \pi_1 (X, x) .
\end{align}

Since the fundamental group of a path-connected space is unique up to isomorphism regardless of the choice of base point, the object groups of a connected groupoid are all isomorphic.  A \emph{deformation retraction}, $r: G \to G(x)$, of a connected groupoid $G$ to its object group $G(x)$ at a chosen base point $x\in P$ is defined by a choice of elements $\a_y \in G(x,y)$ for each $y \in P$, with the understanding that $\a_{x} = \I$.  Then the deformation retraction is defined by the rule that for any $g\in G(y,z)$, $r: g \mapsto \a_z^{-1} g \a_y \in G(x)$.  A deformation retraction thus uses a choice of an isomorphism between each base point to map every element of $G$ to a chosen base point.  The deformation retraction is surjective, so $r(G)=G(x)$ is independent of the choice of elements $\a_y$ used to define it.  We will give an example of a deformation retraction below after we introduce the concept of a graph associated to a groupoid.

A \emph{graph} consists of a set \emph{vertices} and a set of \emph{edges} which are identified with pairs of vertices.   Graphs are \emph{oriented} when all the edges are assigned one of the two possible orientations, i.e., the pairs are considered to be ordered.  It follows that every groupoid defines a graph by assigning a vertex to each object and an oriented edge to each non-identity morphism and simply ignoring the composition between the morphisms.  Typically, such a graph is not very useful because it contains too many edges.  On the other hand, given a subgraph in the graph of a groupoid, there exists a subgroupoid which is \emph{generated} by the subgraph: it is the smallest subgroupoid containing the subgraph, i.e., the intersection of all the subgroupoids which contain the subgraph. One can then look for a smallest possible subgraph which generates the entire groupoid. This motivates looking for a presentation of a groupoid in terms of a generating graph and a set of relations.

A groupoid is said to be \emph{freely generated} by a graph when there are no non-trivial relations between its elements (i.e., beyond the ones imposed by the identity elements and inverses).  Operationally, we construct such a groupoid by considering identities at all the vertices, elevating each edge in the graph to an isomorphism, and allowing for free composition between the elements or their inverses whenever it is well-defined.

The correspondence between graphs and groupoids allows us to introduce the notion of a presentation of a groupoid in terms of a generator graph and a set of relations.  The basic idea is that relations are some given set of words in a groupoid $G$ such that each of these words forms a loop in the graph of the groupoid at some object.   This is described as a wide (containing all the objects) and totally disconnected (there is no edge connecting two distinct vertices) subgraph $R$ of the graph of $G$.  Thus $R := \bigcup_{x\in P} R(x)$, where $R(x)$ are subsets of elements of the object group $G(x)$ for all objects $x$ of $G$, provided that if a given object does not have any such words defined we take the corresponding subset to consist only of the identity element at the object. 

Let $\rho_x$ denote a generic element in $R$, a possibly-trivial loop based at some object $x$.  For each $\rho_x$ we introduce a relation $\rho_x = \I_x$ in $G$.  Introducing these relations in $G$ is equivalent to forming a new groupoid by making the identifications $\eta \rho_x = \eta$ and $\rho_x \zeta = \zeta$ for all $\rho_x\in R$ and for arbitrary elements $\eta$ and $\zeta\in G$ for which the compositions are well defined.  The precise statement is that we form the quotient groupoid $G/N(R)$ where $N(R)$ is the normal closure of $R$; see, e.g., section 8.3 of \cite{Brown06}.  

The object groups of $H=G/N(R)$ are simply related to those of $G$; see, e.g., theorem 8.3.3 of \cite{Brown06}.  Choose any deformation retraction $r: G \to G(x)$, then $H(x)$ is isomorphic to the group $G(x)$ with the relations $r(\rho) = \I$, $\rho \in R$.

For an example of a deformation retraction in a graph, consider the graph
\begin{align}\label{df1}
\begin{tikzpicture}[->,>=stealth',shorten >=1pt,thick,
    auto,node distance=2.0cm,
    main node/.style={circle, fill=black!20},
    baseline={([yshift=-.5ex]current bounding box.center)}]
\tikzset{every loop/.style={min distance=10mm,looseness=10}}
  \node[main node] (1) {$y$};
  \node[main node] (2) [right of=1] {$x$};
  \path
    (2) edge node[above] {$\a_y$} (1)
        edge [in=-30,out=30,loop] node {$\b$} (2)
    (1) edge [in=150,out=210,loop] node {$\g$} (1);
\end{tikzpicture} .
\end{align}
Vertices are gray circles and correspond to objects (base points) and edges correspond to morphisms (elements) of the associated freely-generated groupoid. Note that the inverses correspond to the same edges but with opposite orientation and that identities, which are loops based at a given basepoint, are not explicitly shown.  This graph generates the groupoid
\begin{align}\label{df1.5}
G = \left\langle
\begin{tikzpicture}[->,>=stealth',shorten >=1pt,thick,
    auto,node distance=2.0cm,
    main node/.style={circle, fill=black!20},
    baseline={([yshift=-.5ex]current bounding box.center)}]
\tikzset{every loop/.style={min distance=10mm,looseness=10}}
  \node[main node] (1) {$y$};
  \node[main node] (2) [right of=1] {$x$};
  \path
    (2) edge node[above] {$\a_y$} (1)
        edge [in=-30,out=30,loop] node {$\b$} (2)
    (1) edge [in=150,out=210,loop] node {$\g$} (1);
\end{tikzpicture} 
\ \right\rangle .
\end{align}
This notation, analogous to that used in presentations of groups in terms of generators and relations, means that $G$ is the groupoid with object set $P=\{x,y\}$ and freely generated by the three elements $\a_y$, $\b$, and $\g$.

The deformation retract, $r(G)$, of this groupoid along $\a_y$ is the freely generated groupoid
\begin{align}\label{df2}
r(G) = \left\langle
\begin{tikzpicture}[->,>=stealth',shorten >=1pt,thick,
    auto,node distance=2.0cm,
    main node/.style={circle, fill=black!20},
    baseline={([yshift=-.5ex]current bounding box.center)}]
\tikzset{every loop/.style={min distance=10mm,looseness=10}}
  \node[main node] (6) {$x$};
  \path
    (6) edge [in=-30,out=30,loop] node {$\b$} (6)
        edge [in=150,out=210,loop] node {$\d := \a_y^{-1}\g\a_y$} (6);
\end{tikzpicture} 
\ \right\rangle .
\end{align}
The retracted graph has a single base point, so it generates a group which is the object group corresponding to the free groupoid above, $r(G)=G(x)$.  In this case it is the group freely generated on two generators, $G(x) = \langle \d, \b \rangle \simeq \Z * \Z$.

For an example of the description of a groupoid via a generator graph and a set of relations, consider 
\begin{align}\label{grex}
G = \left\langle
\begin{tikzpicture}[->,>=stealth',shorten >=1pt,thick,
    auto,node distance=2.0cm,
    main node/.style={circle,
    fill=black!20},
    baseline={([yshift=-.5ex]current bounding box.center)}]
\tikzset{every loop/.style={min distance=10mm,looseness=10}}
  \node[main node] (1) {$y$};
  \node[main node] (2) [right of=1] {$x$};
  \path
    (2) edge[bend right] node[above] {$\b$} (1)
        edge [in=30,out=-30,loop] node[right] {$\g$} (2)
    (1) edge[bend right] node[below] {$\a$} (2);
\end{tikzpicture}
\ \Bigg\vert
\ \g^n \ ,\ \g(\a\b)\g^{-1} (\a\b)^{-1} \ 
\right\rangle .
\end{align}
The relations imposed on the free groupoid tell us that $\g$ is an element of order $n$ and that $\g$ and $\a \b$ --- which are elements based at $x$ --- commute.  The object group $G(x)$ of this groupoid is found via a deformation retract, $r$, of the groupoid along $\b$ to be
\begin{align}\label{grpex}
r(G) = \left\langle
\begin{tikzpicture}[->,>=stealth',shorten >=1pt,thick,
    auto,node distance=2.0cm,
    main node/.style={circle, fill=black!20},
    baseline={([yshift=-.5ex]current bounding box.center)}]
\tikzset{every loop/.style={min distance=10mm,looseness=10}}
  \node[main node] (6) {$x$};
  \path
    (6) edge [in=-30,out=30,loop] node {$\g$} (6)
        edge [in=150,out=210,loop] node {$\d := \a \b$} (6);
\end{tikzpicture}
\ \Bigg\vert
\ \g^n \ ,\ \g \d \g^{-1} \d^{-1} \ 
\right\rangle .
\end{align}
Thus $r(G) = G(x) = \vev{\g,\d \,|\, \g^n, \g\d\g^{-1}\d^{-1}} \simeq \Z\oplus\Z_n$ with $\d := \b\a$ generating $\Z$ and $\g$ generating $\Z_n$.  The result is independent of the choice of deformation retract: retracting along $\a^{-1}$ is easily seen to give an equivalent presentation of $G(x)$.

\section{Seifert--van Kampen theorem for groupoids}
\label{sec:maintheorem}

In \cite{Brown06} we find the following generalization of the Seifert--van Kampen theorem formulated in terms of fundamental groupoids.

If a path-connected space $X$ is the union $X=A \cup B$ of two path-connected spaces $A$ and $B$, the groupoid generalization of the Seifert--van Kampen theorem relates the fundamental groupoid of $X$ to those of $A$, $B$, and $C=A\cap B$.   In particular, $C$ need not be path connected.   Choose a finite representative set $P\subset C$ such that the intersection of $P$ with each connected component of $C$ is a single point.  We will choose $P$ as the object set (set of base points) to define the fundamental groupoids $\pi A$, $\pi B$, $\pi C$, and $\pi X$.  The following theorem determines a presentation of the (any) object group of $\pi X$ in terms of presentations of the groupoids $\pi A$, $\pi B$, and $\pi C$.  (In particular, we do not need to reconstruct the full groupoid data of $\pi X$.)

The category-theoretic pushout of set inclusions induces the pushout of fundamental groupoids
\begin{align}
\begin{matrix}
C & \ \ \subset\ \ \ & A \\[1.5mm]
\cap & & \cap\\[1.5mm]
B & \ \ \subset\ \ \ & X 
\end{matrix}
\qquad \to \qquad
\begin{CD}
\pi C @>{i}>> \pi A\\
@V{j}VV  @VV{u}V \\
\pi B @>>{v}> \pi X \\
\end{CD}
\nn
\end{align}
where $i$, $j$, $u$, and $v$ are the morphisms between groupoids induced by the set inclusions.  Note that all the groupoids have the same set $P$ of objects, that $i, j, u, v$ are the identity on objects, and that $\pi A$ and $\pi B$ are connected groupoids and $\pi C$ is totally disconnected. 

With this set up, we have \cite{Brown06}

\paragraph{Theorem 3.1} 
\label{gsvkthm} 

Choose a base point $x_0\in P$, and let $r: A \to A(x_0)$, $s: B \to B(x_0)$ be deformation retractions obtained by choosing elements $\a_y \in A(x_0, y)$, $\b_y \in B(x_0, y)$, for each $y \in P$, with the understanding that $\a_{x_0} = \I_{x_0}, \b_{x_0} = \I_{x_0}$.  Define
\begin{align}\label{fdef}
f_y := (u \a_y)^{-1}(v \b_y) \in \pi X(x_0),
\end{align} 
and let $F(x_0)$ be the free group on the elements $f_y$, for all $y \in P$, with the relation $f_{x_0} = \I$,
\begin{align}\label{Fdef}
F(x_0) := \left\langle f_y , \ \forall y\in P \ \Big\vert \ f_{x_0} \right\rangle .
\end{align}
Then, the object group $\pi X(x_0)$ is isomorphic to the quotient of the free product group $\pi A(x_0) \ast \pi B(x_0) \ast F(x_0) $ by the relations $(ri\g) f_y (sj\g)^{-1} f_y^{-1} = \I$ for all $y \in P$ and all $\g \in \pi C(y)$, which we write as follows
\begin{align}\label{piXx}
\pi X(x_0) = \left\langle \pi A(x_0) \ast \pi B(x_0) \ast F(x_0) \ \Big\vert \ 
(ri\g) f_y (sj\g)^{-1} f_y^{-1} , \ \forall \g\in \pi C(y), \ \forall y\in P \right\rangle .
\end{align}

For the proof we refer the reader to theorem 8.4.1 in \cite{Brown06}.  
We instead illustrate the application of \nameref{gsvkthm} by using it to compute the fundamental group of a 2-torus.  This simple example is meant to help make the theorem's content clear.  

\paragraph{Example 3.2}
\label{Ex-piT2}

\begin{figure}[ht]
\centering
\includegraphics[width=.50\textwidth]{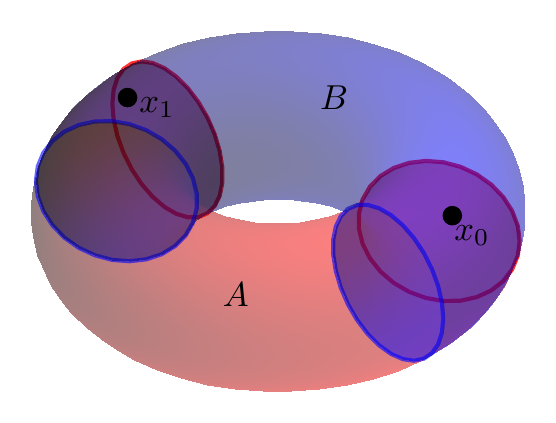}
\caption{A 2-torus, $T$, given as a union of sets $A$ (in red) and $B$ (in blue).  $C=A\cap B$ (in purple) has two components.  A  representative set $P=\{x_0,x_1\}$ of base points is shown.}
\label{torus-fig}
\end{figure}

Realize the 2-torus, $T$, as the union $T=A \cup B$ where we choose the subsets $A$ and $B$ and a representative set of base points $P=\{x_0,x_1\}$ satisfying the conditions of the theorem as in figure \ref{torus-fig}.

With these choices, the fundamental groupoids of the various sets are easily seen to be free, and to be given by
\begin{align} 
\pi A &= \left\langle
\begin{tikzpicture}[->,>=stealth',shorten >=1pt,thick,
    auto,node distance=2.0cm,
    main node/.style={circle, fill=black!20},
    baseline={([yshift=-.5ex]current bounding box.center)}]
\tikzset{every loop/.style={min distance=10mm,looseness=10}}
  \node[main node] (1) {$x_1$};
  \node[main node] (2) [right of=1] {$x_0$};
  \path
    (2) edge node[above] {$\a$} (1)
        edge [in=20,out=-20,loop] node[right] {$\d_A$} (2);
\end{tikzpicture}
\right\rangle , 
\nn\\
\pi B &= \left\langle
\begin{tikzpicture}[->,>=stealth',shorten >=1pt,thick,
    auto,node distance=2.0cm,
    main node/.style={circle,fill=black!20},
    baseline={([yshift=-.5ex]current bounding box.center)}]
\tikzset{every loop/.style={min distance=10mm,looseness=10}}
  \node[main node] (1) {$x_1$};
  \node[main node] (2) [right of=1] {$x_0$};
  \path
    (2) edge node[above] {$\b$} (1)
        edge [in=20,out=-20,loop] node[right] {$\d_B$} (2);
\end{tikzpicture}
\right\rangle , 
\nn\\
\pi C &= \left\langle
\begin{tikzpicture}[->,>=stealth',shorten >=1pt,thick,
    auto,node distance=2.0cm,
    main node/.style={circle,fill=black!20},
    baseline={([yshift=-.5ex]current bounding box.center)}]
\tikzset{every loop/.style={min distance=10mm,looseness=10}}
  \node[main node] (1) {$x_1$};
  \node[main node] (2) [right of=1] {$x_0$};
  \path
    (1) edge [in=160,out=200,loop] node[left] {$\g_1$} (1)
    (2) edge [in=20,out=-20,loop] node[right] {$\g_0$} (2);
\end{tikzpicture}
\right\rangle . \nn
\end{align}
Furthermore, their object groups are
\begin{align} 
\pi A(x_0) &= \vev{\d_A} \simeq \Z, &
\pi B(x_0) &= \vev{\d_B} \simeq \Z , \nn\\
\pi C(x_0) &= \vev{\g_0} \simeq \Z , &
\pi C(x_1) &= \vev{\g_1} \simeq \Z . \nn
\end{align}

The inclusion maps $C\subset A$ and $C\subset B$ induce the morphisms between the associated fundamental groupoids which act on the generating elements of $\pi C$ by
\begin{align} \label{T2ij}
i(\g_0) &= \d_A, &    j(\g_0) &= \d_B , \nn\\
i(\g_1) &= \a\d_A\a^{-1}, &    j(\g_1) &=  \b\d_B\b^{-1} . 
\end{align}
Define deformation retractions $r:\pi A\to\pi A(x_0)$ and $s:\pi B\to\pi B(x_0)$ by choosing elements $\{ \a_{x_0} = \I_{x_0}, \a_{x_1} = \a \} \subset \pi A$ to define $r$ and $\{ \b_{x_0}=\I_{x_0}, \b_{x_1} = \b \} \subset \pi B$ to define $s$.  It then follows, using \eqref{T2ij} and similar to \eqref{df2}, that
\begin{align} \label{T2risj}
ri(\g_0) &= \d_A ,&
sj(\g_0) &= \d_B, \nn\\
ri(\g_1) &= \a^{-1} i(\g_1)\a = \d_A , &
sj(\g_1) &= \b^{-1} j(\g_1)\b = \d_B . 
\end{align}

Define, using the prescription \eqref{Fdef}, the group $F(x_0)$ generated by $f_y$ given in \eqref{fdef} for all $y\in P$.  This gives $F(x_0) = \vev{f_0,\ f_1\ |\ f_0}$ where $f_0 := f_{x_0} = (u \I_{x_0})^{-1}(v \I_{x_0})$ and $f_1 := f_{x_1} = (u\a)^{-1} (v\b)$.  Since $u$ and $v$ are morphisms, $f_0$ must be the identity element at $x_0$ in $\pi T$.  Indeed, the relation $f_0$ identifies it with the group identity, so
\begin{align} 
F(x_0) &\simeq \vev{f_1} \simeq\Z. \nn
\end{align}
Since a representative of $\a$ is a path from $x_0$ to $x_1$ in $A$ and similarly $\b$ is represented by a path from $x_0$ to $x_1$ in $B$, 
the ``new'' generator $f_1$ introduced by the prescription of \nameref{gsvkthm} has as a representative a closed path winding once around the longitudinal cycle of $T$.  Constructing such representative paths for the $f_y$ generators is not required to carry out the algebraic construction of $\pi T(x_0)$ given in the theorem;  we have given it here only to help those readers who would like to visualize the construction.

Now apply \eqref{piXx} to compute the fundamental group of the 2-torus.  The main work is to compute the relations for all $y\in P = \{x_0,x_1\}$ and all $\g\in\pi C(y)$.  For $y=x_0$, $\g\in\pi C(y)$ are the elements $\g = (\g_0)^n$ for $n\in\Z$, so the associated relations are
\begin{align} 
(ri\g_0^n)f_0 (sj\g_0^n)^{-1} f_0^{-1}
= (ri\g_0)^n f_0 (sj\g_0)^{-n} f_0^{-1}, 
\qquad n\in\Z, \nn
\end{align}
where the equality follows since $i$,$j$,$r$, and $s$ act as homomorphisms on object groups.  In particular, this implies that the members of this infinite set of relations are not independent since they all follow from the relation with $n=1$.  Using the explicit form of $ri$ and $sj$ calculated in \eqref{T2risj} then gives the single $y=x_0$ relation $(ri\g_0)f_0 (sj\g_0)^{-1} f_0^{-1} = \d_A f_0 \d_B^{-1} f_0^{-1}$.  A similar computation for $y=x_1$ gives one new relation $\d_A f_1 \d_B^{-1} f_1^{-1}$.  Thus by \nameref{gsvkthm}
\begin{align} 
\pi_1(X) &= \pi X(x_0) = 
\vev{\pi A(x_0)\ast \pi B(x_0) \ast F(x_0)\  
|\  \d_A f_0 \d_B^{-1} f_0^{-1},
\ \d_A f_1 \d_B^{-1} f_1^{-1}      }
\nn\\
&= 
\vev{\d_A,\ \d_B,\ f_0,\ f_1\ |\ f_0,\ \d_A f_0 \d_B^{-1} f_0^{-1},
\ \d_A f_1 \d_B^{-1} f_1^{-1} }
\nn\\
&= 
\vev{\d_A,\ \d_B,\ f_1\ |\ \d_A \d_B^{-1},
\ \d_A f_1 \d_B^{-1} f_1^{-1} }
= 
\vev{\d_A,\ f_1\ |\ \ \d_A f_1 \d_A^{-1} f_1^{-1} }
\simeq \Z\oplus\Z .
\nn
\end{align}
This is the well known result for the fundamental group of the 2-torus.

Before moving on to the main application to torus links, some remarks may be useful.  (1) The relations in \eqref{piXx} with $y=x_0\in P$ may be non-trivial and thus must be kept, even though $f_{x_0} = \I$ in \eqref{Fdef}.  (2) The relations $(ri\g) f_y (sj\g)^{-1} f_y^{-1} = \I$ are trivial for $\g=\I\in C(y)$ for any $y$, so we need not include them.  (3) More generally, it is sufficient to include only the relations coming from the $\g$ in a generating subset of $\pi C(y)$.  (4) The number of independent relations can be reduced by choosing the deformation retractions $r$ and $s$ intelligently with respect to a given presentation of $\pi A$ and $\pi B$.

\section{Knot group of torus links}
\label{sec:links}

Consider a torus link $\kn$ as defined in \eqref{Knqpdef}, and sets $A$ and $B$ defined in \eqref{Knpqsets} whose union, $X$, is the complement of $\kn$ in $S^3$.   The intersection $C=A\cap B$ has $n$ disjoint path components.  We select one base point in each of those path-connected subspaces and index them using integers modulo $n$ to form the representative object set $P : = \{x_0, \dots, x_{n-1} \}$ on which to define all the fundamental groupoids.   It will be convenient to make a concrete choice of the $x_k$ in such a way as to make the deformation retractions used below have simple forms.   One such choice is to take
\begin{align}\label{specificP}
x_k 
:= \left\{ (\chi,\th,\psi) 
= ( \tfrac\pi4,\ (2k+1)\tfrac{\pi}{nq},\ k \tfrac{\pi}{np} ) \right\}, 
\qquad  k\in\{0,\ldots,n{-}1\},
\end{align}
in the coordinates defined in section \ref{sec:definitions}.  It is easy to check that with this definition, each $x_k$ lies in a different path component of $C$.

Since each path component of $C$ deformation retracts to $\kone$ which is homeomorphic to $S^1$, the object group $\pi C(x) \simeq \Z$ for any $x\in P$.  Since $\pi C$ is totally disconnected, it is freely generated by $n$ homotopy cycles, one at each base point.  The spaces $A$ and $B$ have isomorphic fundamental groupoids. These groupoids are generated by a set of elements that connect base points $x_k$ to $x_{k+1}$ with the subscripts identified modulo $n$.  This is stated more precisely as

\paragraph{Lemma 4.1} 
\label{linkgraph} 

The following graphs freely generate the fundamental groupoids $\pi A$, $\pi B$, and $\pi C$:
\begin{align} 
\pi A &= \left\langle
\scalebox{0.75}{%
\begin{tikzpicture}[->,>=stealth',shorten >=1pt,thick,
    auto,node distance=2.0cm,
    main node/.style={circle, fill=black!20},
    baseline={([yshift=-.5ex]current bounding box.center)}]
\tikzset{every loop/.style={min distance=10mm,looseness=10}}
  \node[main node] (0) at (0,3) {$x_0$};
  \node[main node] (1) at (1.5,2.1) {$x_1$};
  \node (2) at (.5,.9) {};
  \node at (.35,.85) {$\bullet$};
  \node at (0,.8) {$\bullet$};
  \node at (-.35,.85) {$\bullet$};
  \node (n-2) at (-.5,.9) {};
  \node[main node] (n-1) at (-1.5,2.1) {$\phantom{x_0}$};
  \node at (-1.5,2.1) {\footnotesize $x_{n{-}1}$};
  \path
    (0) edge[bend left] node[above] {$\ta_1$} (1)
    (1) edge[bend left] node[above left] {$\ta_2$} (2)
    (n-2) edge[bend left] node[above right] {$\ta_{n{-}1}$} (n-1)
    (n-1) edge[bend left] node[above] {$\ta_n$} (0);
\end{tikzpicture} }
\right\rangle ,& &&
\pi B &= \left\langle
\scalebox{0.75}{%
\begin{tikzpicture}[->,>=stealth',shorten >=1pt,thick,
    auto,node distance=2.0cm,
    main node/.style={circle, fill=black!20},
    baseline={([yshift=-.5ex]current bounding box.center)}]
\tikzset{every loop/.style={min distance=10mm,looseness=10}}
  \node[main node] (0) at (0,3) {$x_0$};
  \node[main node] (1) at (1.5,2.1) {$x_1$};
  \node (2) at (.5,.9) {};
  \node at (.35,.85) {$\bullet$};
  \node at (0,.8) {$\bullet$};
  \node at (-.35,.85) {$\bullet$};
  \node (n-2) at (-.5,.9) {};
  \node[main node] (n-1) at (-1.5,2.1) {$\phantom{x_0}$};
  \node at (-1.5,2.1) {\footnotesize $x_{n{-}1}$};
  \path
    (0) edge[bend left] node[above] {$\tb_1$} (1)
    (1) edge[bend left] node[above left] {$\tb_2$} (2)
    (n-2) edge[bend left] node[above right] {$\tb_{n{-}1}$} (n-1)
    (n-1) edge[bend left] node[above] {$\tb_n$} (0);
\end{tikzpicture} }
\right\rangle ,
\nn\\
\pi C &= \left\langle
\scalebox{0.75}{%
\begin{tikzpicture}[->,>=stealth',shorten >=1pt,thick,
    auto,node distance=2.0cm,
    main node/.style={circle, fill=black!20},
    baseline={([yshift=-.5ex]current bounding box.center)}]
\tikzset{every loop/.style={min distance=10mm,looseness=10}}
  \node[main node] (0) at (0,3) {$x_0$};
  \node[main node] (1) at (1.5,2.1) {$x_1$};
  \node (2) at (.5,.9) {};
  \node at (.7,1.1) {$\bullet$};
  \node at (.25,.85) {$\bullet$};
  \node at (-.25,.85) {$\bullet$};
  \node at (-.7,1.1) {$\bullet$};
  \node (n-2) at (-.5,.9) {};
  \node[main node] (n-1) at (-1.5,2.1) {$\phantom{x_0}$};
  \node at (-1.5,2.1) {\footnotesize $x_{n{-}1}$};
  \path
    (0) edge[loop above] node[above] {$\g_0$} (0)
    (1) edge[in=30,out=60,loop] node[above left] {$\g_1$} (1)
    (n-1) edge[in=120,out=150,loop] node[above] {$\g_{n-1}$} (n-1);
\end{tikzpicture} }
\right\rangle .
\nn
\end{align} 
\noindent \textbf{Proof:}  $A$ is homotopic to a circle $S^1$ with $n$ base points.  Concretely, the deformation retraction which retracts the $0\le \chi \le \pi/4$ interval to $\chi=0$ while acting as the identity on the $\th,\ps$ coordinates, retracts $A$ to the ``interior'' circle $S^1_A$ defined in section \ref{sec:definitions}.  Under this retraction the base points map to the points $x_k \mapsto (\chi=0,\, \th=(2k+1)\tfrac{\pi}{nq})$ on the circle.  In particular, with the choice \eqref{specificP} they are ordered along the circle with increasing $k$ (mod $n$).  Define $\ta_k$ to be the homotopy equivalence class of the path on $S^1_A$ from $x_{k-1}$ to $x_k$ consisting of the interval $(2k-1)\frac\pi{nq} \le \th \le (2k+1)\frac\pi{nq}$ for $k=1,\ldots, n-1$, and the interval $(2n-1)\frac\pi{nq}\le \th \le \frac\pi{nq}$ for $k=n$ (where $\th$ is periodic with period $2\pi$).  The graph associated to these elements is shown above.  An arbitrary element of $\pi A$ connects some $x_i$ to $x_j$ by winding $m \in \Z$ number of times around the circle.  More precisely, this means that this element can be written as a product of the $\ta_j$ elements of the form $\ta_{j}  \dots  \ta_{i+2} \ta_{i+1}  (\ta_i \dots \ta_{i+2} \ta_{i+1})^m$. The set of such words is in one to one correspondence with the set of elements of the groupoid freely generated by the graph.   The argument for $B$ is similar, where one now retracts to the ``exterior'' $S^1_B$ circle by mapping $\pi/4\le\chi\le\pi/2$ to $\chi=\pi/2$ without affecting $\th$ and $\ps$.  With the choice \eqref{specificP}, the base points retract to points ordered along $S^1_B$ with increasing $j$, and so gives the generating graph shown in a manner completely analogous to that of $\pi A$.  The case of $C$ is simpler, and was outlined above.  \hfill $\Box$

The concrete choice of base points and of retractions to $S^1_A$ and $S^1_B$ made in the proof of \nameref{linkgraph} was not essential.  Any other choice (which retracts the base points to distinct points) gives the same result, though perhaps with the ordering of the base points along the circles permuted.  It is merely convenient to retract in such a way that the base points are both ordered by increasing $k$ along both $S^1_A$ and $S^1_B$.

Given these groupoids, we fix the base point $x_0$ as in \nameref{gsvkthm}, and define elements $\a_y$ and $\b_y$ for every object $y$ as follows.  We agree to have $\a_{x_0} = \I_{x_0}$ and $\b_{x_0} = \I_{x_0}$. For any object $x_k$, with $k \neq 0$, let $\a_{x_k} : = \ta_k \dots \ta_2 \ta_1$ and $\b_{x_k} : = \tb_k \dots \tb_2 \tb_1$. We carry out the associated retractions $r: \pi A \to \pi A(x_0)$ and $s: \pi B \to \pi B(x_0)$ to get the following object groups.

\paragraph{Lemma 4.2}
\label{linkretract} 

The object groups of $\pi A$ and $\pi B$ at $x_0$ are $\pi A(x_0) = \vev\a \cong \Z $ and $\pi B(x_0) = \vev\b  \cong \Z$.

\noindent \textbf{Proof:} The objects group $\pi A(x_0)$ is the subgroupoid of the free groupoid $\pi A$ consisting of all the elements based at $x_0$. From the proof of \nameref{linkgraph} these are the elements of the form $(\ta_n \ta_{n-1} \dots \ta_2 \ta_1)^m$. So $\pi A(x_0)$ is generated by $\a := \ta_n \dots \ta_1$. Similarly,  $\pi B(x_0)$ is generated by $\b := \tb_n \dots \tb_1$. \hfill $\Box$

In the proof of \nameref{linkretract} we introduced the generators $\a$ and $\b$ of $\pi A(x_0)$ and $\pi B(x_0)$.  It will be useful for our considerations in section \ref{sec:nested} to locate homotopic representative paths for them in $A$ and $B$.   As noted in the proof of \nameref{linkgraph} $A$ retracts to the circle $S^1_A$ so $\a=\ta_n \dots \ta_1$ is the class of the path that winds around that circle once.  We can continuously deform this circular path, which is at $\chi=0$ for $\th\in[0,2\pi)$, to a circular path near the boundary, $T$, of $A$ by continuously changing $\chi$ to $\chi = \pi/4 - \epsilon$ for an arbitrarily small positive number $\epsilon$, while keeping $\ps=0$ fixed.  Thus
\begin{align}\label{alpharep}
\a \simeq \{ \chi = \tfrac\pi4-\epsilon, \ 0\le\th<2\pi,\ \ps=0 \}.
\end{align}
Likewise, $\b$ is homotopic to $S^1_B$ which can be continuously deformed to the circle 
\begin{align}\label{betarep}
\b \simeq \{ \chi = \tfrac\pi4+\epsilon, \ \th=\tfrac\pi{nq}, \ 0\le\ps<2\pi\}.
\end{align}
The resulting representative paths for $\a$ and $\b$ are illustrated in figure \ref{fig1}. \hfill $\Box$

\begin{figure}[ht]
\centering
\includegraphics[width=.75\textwidth]{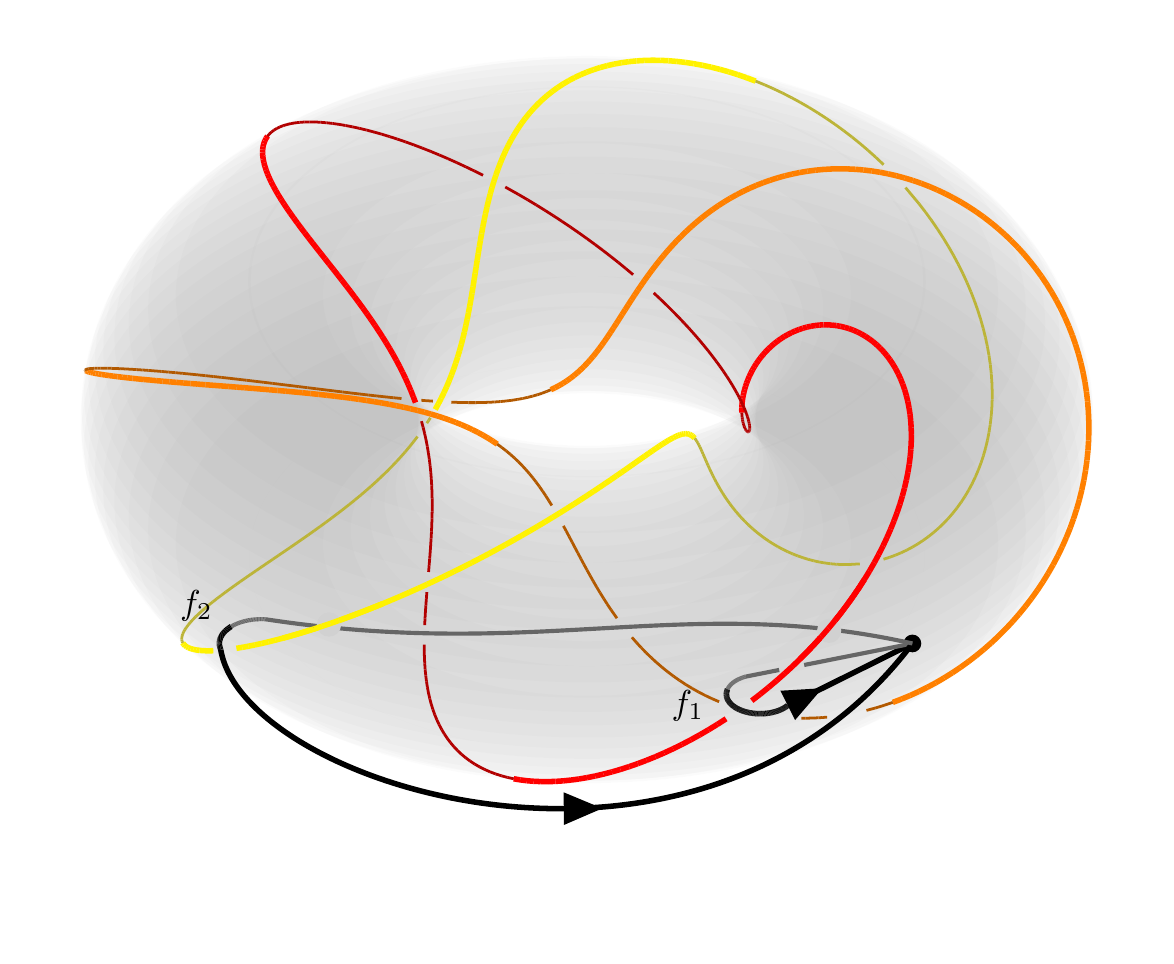}
\caption{A torus 3-link (its components are shown in red, orange, and yellow) with the non-trivial $f_i$ homotopy cycles in $\pi X(x_0)$, where $X$ is the complement in $S^3$ of the torus 3-link.  The shaded solid donut is for visualization purposes only.}
\label{fig3}
\end{figure}

The group, $F(x_0)$, containing the additional generators in the object group $\pi X(x_0)$ required by the pushout, is defined in \nameref{gsvkthm} to be the group generated by the collection of elements $f_k$ --- with the one relation $f_0 = \I$ and free otherwise --- where $f_k := f_{x_k} = (u \a_{x_k})^{-1}(v \b_{x_k})$ for each object $x_k\in P$. That is, $F(x_0) = \vev{ f_0 \, , \, f_1 \, , \, \dots \, , \, f_{n-1} \, | \, f_0 = \I \,} \cong \vev{ f_1 \, , \, \dots \, , \, f_{n-1}}$.  Figure \ref{fig3} illustrates representative paths corresponding to these generators:  they are paths which link up to $n-1$ of strands corresponding to distinct copies of $\kone$ in the $n$-link.  But we emphasize that constructing such representative paths is not needed to carry out the construction of $\pi X(x_0)$ that follows.

The free product of the object groups, $\pi A(x_0) \ast \pi B(x_0) \ast F(x_0)$, is simply presented by forming the unions of all their respective generators and relations.  The fundamental group is formed from this group by adding the relations between the elements of distinct factors in the free product induced by the pushout.  Specifically, \nameref{gsvkthm} instructs us to construct relations of the form $(ri\g) f_k (sj\g)^{-1} f_k^{-1} = \I$, for all base points $x_k\in P$ and for a generating set of elements $\g\in\pi C(x_k)$.  By \nameref{linkgraph} there is only one generator $\g_k\in \pi C(x_k)$ for each $k$ and, furthermore, $\g_k$ deformation retracts to a homotopic copy of $\kone$ in the knot complement.   With the choice of base points given in \eqref{specificP}, representatives of the $\g_k$ homotopy classes are the paths $\g_k \simeq \{ (\chi, \th, \ps) = x_k + (0,\, p\tau,\, q\tau) \ \text{for}\ 0\le\tau<2\pi \}$.  The image of $\g_k$ under the deformation retraction of $A\to S^1_A$ used in \nameref{linkgraph} is $\g_k \simeq \{ (\chi, \th) = (0,\, (2k+1)\frac\pi{nq} + p\tau) \ \text{for}\ 0\le\tau<2\pi \}$, which is a path which starts at (the deformation retract image of) $x_k$ and winds $q$ times around $S^1_A$ in the direction of increasing $\th$.  The inclusion map $i: \pi C \to \pi A$ on the groupoids thus maps $\g_k \mapsto i \g_k = (\ta_{k-1} \dots \ta_{k+1} \ta_k)^p$, where the subscripts are understood to be in the remainder classes modulo $n$.  A similar construction under the inclusion map $j: \pi C \to \pi B$ for the exterior space gives us $\g_k \mapsto j \g_k = (\tb_{k-1} \dots \tb_{k+1} \tb_k)^q$. 

This presentation of $i \g_k$ and $j \g_k$ as a product of elements of $\pi A$ and $\pi B$ is particularly convenient because the retractions defined above were constructed along a similar product of elements.  In particular,  it immediately follows that under the retractions $i \g_k \mapsto r i \g_k = \a^p$ and $j \g_k \mapsto s j \g_k = \b^q$ for all $k\in\Z_n$.  Therefore the relations are
\begin{align}\label{conjugacyrelations}
(ri\g) f_k (sj\g)^{-1} f_k^{-1} & = \a^p f_k \b^{-q} f_k^{-1} 
  \qquad \forall k \in \Z_n.
\end{align}
It now follows from equation \eqref{piXx} that a presentation of the fundamental group of the link complement, based at point $x_0$, is

\paragraph{Theorem 4.3}
\label{result1} 

$\pi X(x_0) = \pi_1 (S^3 \setminus \kn)
=\Bigl\langle \a \, , \b \, , f_k \, \Big| \, 
f_0 \, , \a^p f_ k \b^{-q} f_k^{-1} 
\quad \forall k\in\Z_n  \, \Bigr\rangle$.

\vspace{3mm}
As a consistency check, we note that for $n = 1$ the fundamental group reduces to the torus knot result \eqref{centuryold}.

\section{Knot group of torus links with two unknots}
\label{sec:withunknots}

\begin{figure}[ht]
\centering
\includegraphics[width=.75\textwidth]{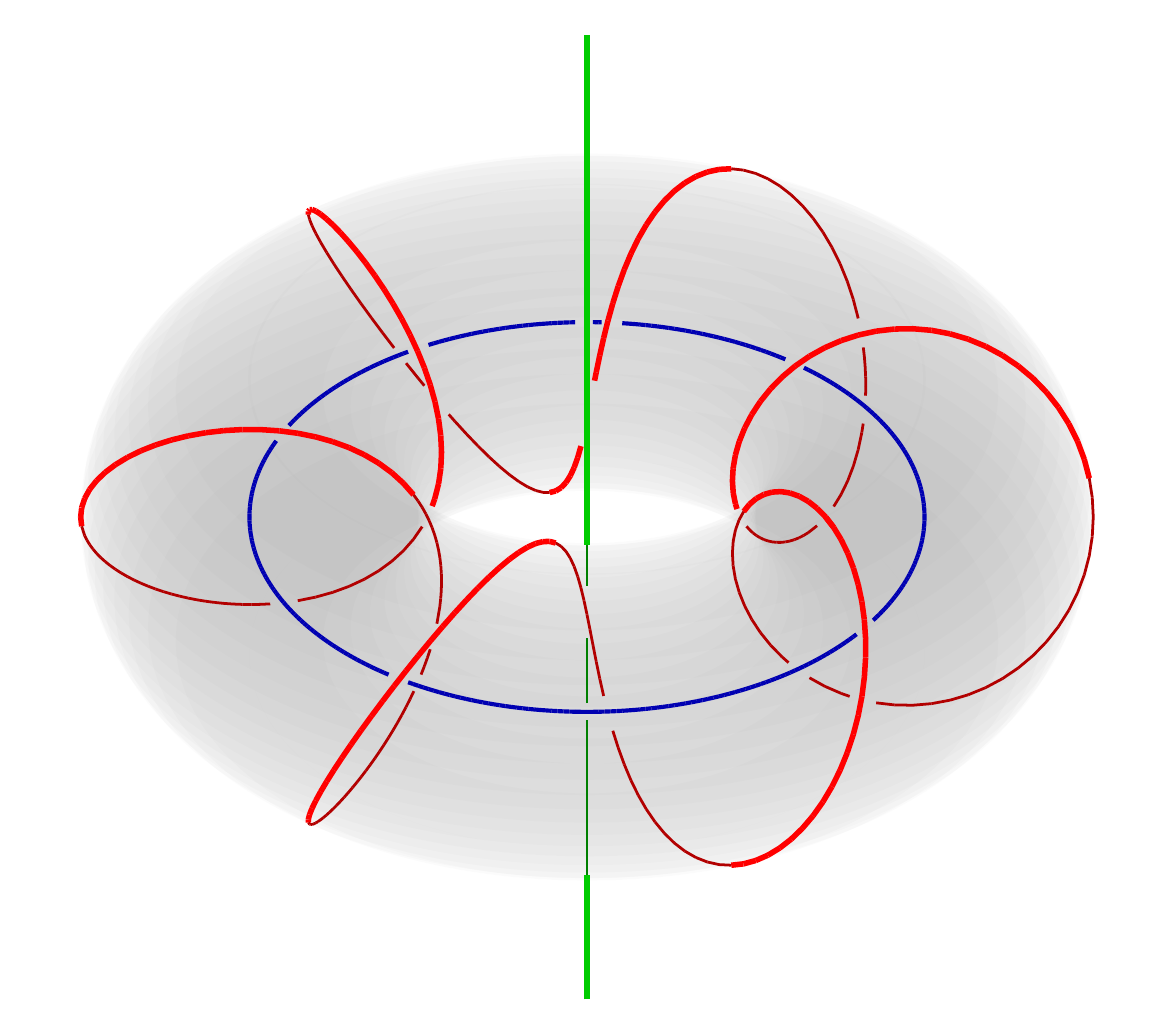}
\caption{A torus knot (red) linked with two unknots (green and blue).  The interior unknot ($\simeq S^1_A$) is blue, and the exterior unknot ($\simeq S^1_B$) is green.  The green unknot looks like a line instead of a circle because it goes through the point on $S^3$ that lies at infinity in the stereographic projection to $\R^3$.  The particular link shown is denoted $S^1_B\cup K^1_{1,6} \cup S^1_A$.}
\label{fig4}
\end{figure}

We now extend the analysis to a torus link further linked with two unknots: one passing through the interior of the auxiliary torus and another one through the exterior.  An arrangement of a knot linked with an interior and an exterior unknot is illustrated in figure \ref{fig4}.  In terms of the explicit parameterization introduced at the beginning of section \ref{sec:definitions}, the interior unknot is the circle $S^1_A$ and the exterior unknot is the circle $S^1_B$, so a torus link linked with these two unknots is $S^1_B \cup \kn \cup S^1_A$.

The computation of the knot group for this case is similar to the one in section \ref{sec:links}.  The consequence of having the two additional unknots on the fundamental groupoids of subsets $A$, $B$, and $C=A\cap B$, is summarized in \nameref{lemma5.2}.  As a preliminary step, we compute the fundamental groupoid of a 2-torus with $n$ base points.

\paragraph{Lemma 5.1}
\label{torusgroupoid} 

The fundamental groupoid for torus $T$ with $n$ base points is given by the following generator graph
\begin{align} 
\pi T &= \Bigg\langle 
\scalebox{0.75}{%
\begin{tikzpicture}[->,>=stealth',shorten >=1pt,thick,
    main node/.style={circle,fill=black!20},
    baseline={([yshift=-.5ex]current bounding box.center)}]
\tikzset{every loop/.style={min distance=10mm,looseness=10}}
  \node[main node] (0) at (0,3) {$x_0$};
  \node[main node] (1) at (1.5,2.1) {$x_1$};
  \node (2) at (.5,.9) {};
  \node at (.35,.85) {$\bullet$};
  \node at (0,.8) {$\bullet$};
  \node at (-.35,.85) {$\bullet$};
  \node (n-2) at (-.5,.9) {};
  \node[main node] (n-1) at (-1.5,2.1) {$\phantom{x_0}$};
  \node at (-1.5,2.1) {\footnotesize $x_{n{-}1}$};
  \path
    (0) edge[loop above] node[above] {$\d$} (0)
    (0) edge[bend left] node[above] {$\ta_1$} (1)
    (1) edge[bend left] node[above left] {$\ta_2$} (2)
    (n-2) edge[bend left] node[above right] {$\ta_{n{-}1}$} (n-1)
    (n-1) edge[bend left] node[above] {$\ta_n$} (0);
\end{tikzpicture} }
\ \Bigg| \ 
(\ta_n\cdots\ta_2\ta_1) \d (\ta_n\cdots\ta_2\ta_1)^{-1} \d^{-1}
\ \Bigg\rangle .\nn
\end{align}
\noindent \textbf{Proof:} Let $x_k$, $k\in\Z_n$, label the $n$ base points in the object set of $\pi T$.  Since $T$ is connected, all the object groups in $\pi T$ are isomorphic to $\pi T(x_0) \simeq \Z\oplus\Z$, which is generated by two commuting generators.  In particular, $\pi T(x_0) \simeq \vev{\a,\d\,|\,\a\d\a^{-1}\d^{-1}}$.  The groupoid elements in $\pi T(x_i,x_j)$ are equivalence classes of paths connecting $x_i$ to $x_j$ and winding $(p,q)$ times around $T$.  Specifically, pick elements $\ta_k \in \pi T(x_k,x_{k+1})$ for each $k\in\Z_n$.  Then the elements of $\pi T(x_i,x_j)$ are $\ta_{j-1}\cdots\ta_1 \d^q \a^p \ta_n \cdots \ta_i$ where we identify $\a = \ta_n \cdots \ta_1$.  These are all the inequivalent words (elements) in the presentation given in the statement of the lemma. \hfill $\Box$

To compute the knot group, we choose, just as before, $A$ and $B$ to be the intersections of the knot complement with the closed ``interior'' and closed ``exterior'', respectively, of the 2-torus $T$ on which the torus knot component is placed.  The (topological) interior of $A$ is thus the ``interior'' of $T$ minus the circle $S^1_A$, while the interior of $B$ is the ``exterior'' of $T$ minus the circle $S^1_B$.

\paragraph{Lemma 5.2} 
\label{lemma5.2}

The fundamental groupoids $\pi A$, $\pi B$, and $\pi C$ have the presentations
\begin{align}  
\pi A &= \Bigg\langle
\scalebox{0.75}{%
\begin{tikzpicture}[->,>=stealth',shorten >=1pt,thick,
    auto,node distance=2.0cm,
    main node/.style={circle, fill=black!20},
    baseline={([yshift=-.5ex]current bounding box.center)}]
\tikzset{every loop/.style={min distance=10mm,looseness=10}}
  \node[main node] (0) at (0,3) {$x_0$};
  \node[main node] (1) at (1.5,2.1) {$x_1$};
  \node (2) at (.5,.9) {};
  \node at (.35,.85) {$\bullet$};
  \node at (0,.8) {$\bullet$};
  \node at (-.35,.85) {$\bullet$};
  \node (n-2) at (-.5,.9) {};
  \node[main node] (n-1) at (-1.5,2.1) {$\phantom{x_0}$};
  \node at (-1.5,2.1) {\footnotesize $x_{n{-}1}$};
  \path
    (0) edge[loop above] node[above] {$\d_A$} (0)
    (0) edge[bend left] node[above] {$\ta_1$} (1)
    (1) edge[bend left] node[above left] {$\ta_2$} (2)
    (n-2) edge[bend left] node[above right] {$\ta_{n{-}1}$} (n-1)
    (n-1) edge[bend left] node[above] {$\ta_n$} (0);
\end{tikzpicture} }
\ \Bigg| \ 
(\ta_n\cdots\ta_2\ta_1) \d_A (\ta_n\cdots\ta_2\ta_1)^{-1} \d_A^{-1}
\ \Bigg\rangle,
\nn\\
\pi B &= \Bigg\langle
\scalebox{0.75}{%
\begin{tikzpicture}[->,>=stealth',shorten >=1pt,thick,
    auto,node distance=2.0cm,
    main node/.style={circle,
    fill=black!20},
    baseline={([yshift=-.5ex]current bounding box.center)}]
\tikzset{every loop/.style={min distance=10mm,looseness=10}}
  \node[main node] (0) at (0,3) {$x_0$};
  \node[main node] (1) at (1.5,2.1) {$x_1$};
  \node (2) at (.5,.9) {};
  \node at (.35,.85) {$\bullet$};
  \node at (0,.8) {$\bullet$};
  \node at (-.35,.85) {$\bullet$};
  \node (n-2) at (-.5,.9) {};
  \node[main node] (n-1) at (-1.5,2.1) {$\phantom{x_0}$};
  \node at (-1.5,2.1) {\footnotesize $x_{n{-}1}$};
  \path
    (0) edge[loop above] node[above] {$\d_B$} (0)
    (0) edge[bend left] node[above] {$\tb_1$} (1)
    (1) edge[bend left] node[above left] {$\tb_2$} (2)
    (n-2) edge[bend left] node[above right] {$\tb_{n{-}1}$} (n-1)
    (n-1) edge[bend left] node[above] {$\tb_n$} (0);
\end{tikzpicture} }
\ \Bigg| \ 
(\tb_n\cdots\tb_2\tb_1) \d_B (\tb_n\cdots\tb_2\tb_1)^{-1} \d_B^{-1}
\ \Bigg\rangle ,
\nn\\
\pi C &= \left\langle
\scalebox{0.75}{%
\begin{tikzpicture}[->,>=stealth',shorten >=1pt,thick,
    auto,node distance=2.0cm,
    main node/.style={circle,
    fill=black!20},
    baseline={([yshift=-.5ex]current bounding box.center)}]
\tikzset{every loop/.style={min distance=10mm,looseness=10}}
  \node[main node] (0) at (0,3) {$x_0$};
  \node[main node] (1) at (1.5,2.1) {$x_1$};
  \node (2) at (.5,.9) {};
  \node at (.7,1.1) {$\bullet$};
  \node at (.25,.85) {$\bullet$};
  \node at (-.25,.85) {$\bullet$};
  \node at (-.7,1.1) {$\bullet$};
  \node (n-2) at (-.5,.9) {};
  \node[main node] (n-1) at (-1.5,2.1) {$\phantom{x_0}$};
  \node at (-1.5,2.1) {\footnotesize $x_{n{-}1}$};
  \path
    (0) edge[loop above] node[above] {$\g_0$} (0)
    (1) edge[in=30,out=60,loop] node[above left] {$\g_1$} (1)
    (n-1) edge[in=120,out=150,loop] node[above] {$\g_{n-1}$} (n-1);
\end{tikzpicture} }
\right\rangle .\nn
\end{align} 
\noindent \textbf{Proof:}  $A$ and $B$ are each homotopic to a 2-torus.  For example, for $A$ a homotopy is given by retracting the coordinate interval $0\le\chi < \pi/4$ to the value $\chi=\pi/4$ in the parameterization given in the beginning of section \ref{sec:definitions}, and similarly for $B$.  Thus, their fundamental groupoids with $n$ base points is as described in \nameref{torusgroupoid}.  The intersection fundamental groupoid for $C$ is the same as in \nameref{linkgraph}. \hfill $\Box$

Once again, we define retractions $r: \pi A \to \pi A(x_0)$ and $s: \pi B \to \pi B(x_0)$ using the elements $\a_y$ and $\b_y$ defined for every object $y$ as follows: for $y = x_0$, $\a_{x_0} = \I_{x_0}$ and $\b_{x_0} = \I_{x_0}$; for $y = x_k \neq x_0$, $\a_{x_k} : = \ta_k \dots \ta_2 \ta_1$ and $\b_{x_k} : = \tb_k \dots \tb_2 \tb_1$. 

\paragraph{Lemma 5.3}
\label{unknotretract}

The object groups at $x_0$ of $\pi A$ and $\pi B$ are given by
%
$\pi A(x_0) = \langle \a \, , \d_A \, | $ $\a \d_A \a^{-1} \d_A^{-1} \, \rangle $ and $\pi B(x_0) = \langle \b \, , \d_B \, | \b \d_B \b^{-1} \d_B^{-1} \, \rangle $.

\noindent \textbf{Proof:}  This follows from choosing retractions as in the proof of \nameref{linkretract}.  In particular, since the spaces $A$ and $B$ both deformation retract to a torus $T$, the object groups are isomorphic to the one we obtained in \nameref{Ex-piT2}.  \hfill $\Box$

We can summarize the above lemmas by saying that the existence of additional unknots introduces in the corresponding object groups additional non-trivial loops --- the generators $\d_A$ and $\d_B$ --- that link these unknots, with the property that the generators commute with the other generators of the respective object groups.  

\begin{figure}[ht]
\centering
\includegraphics[width=.75\textwidth]{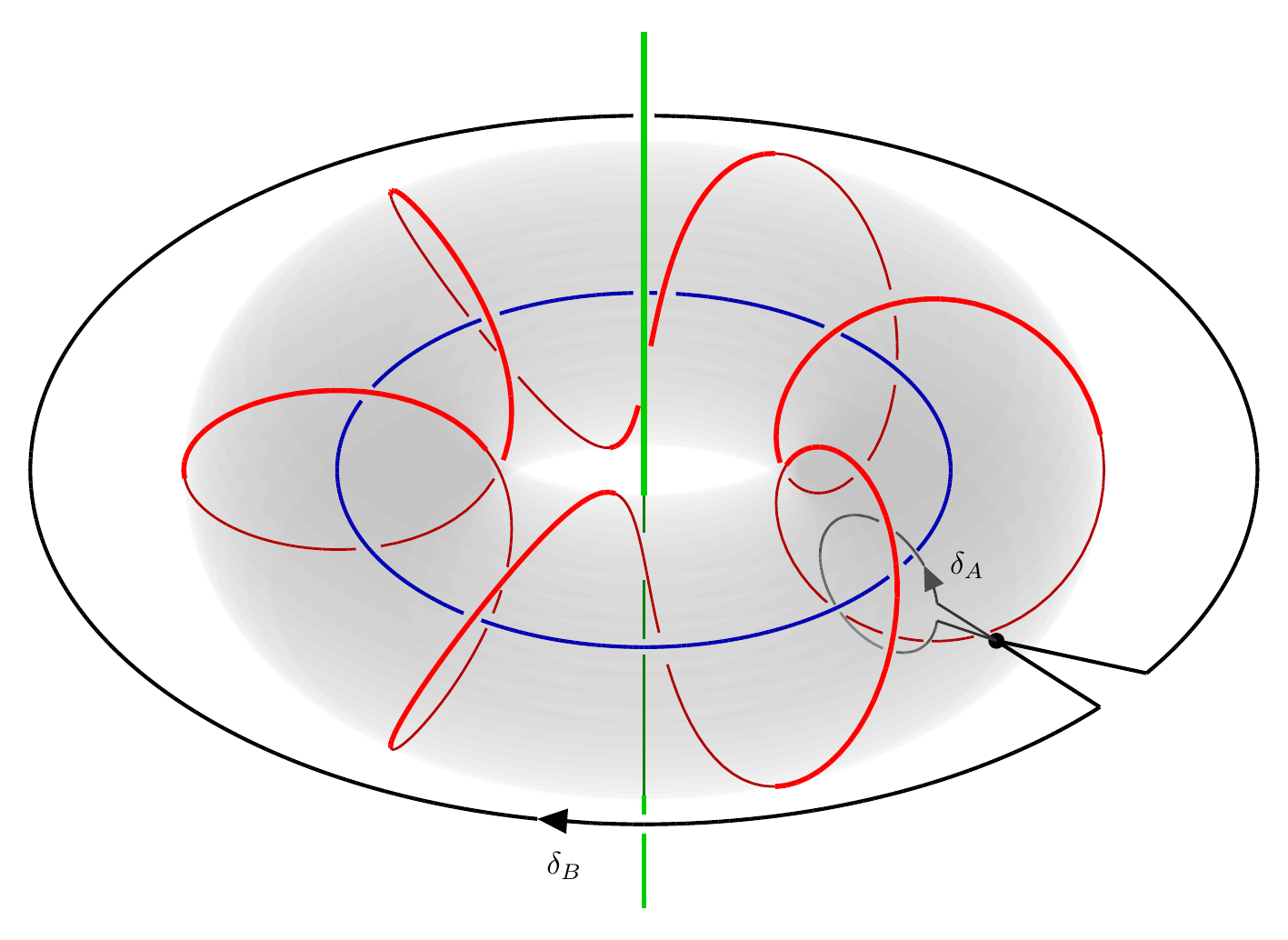}
\caption{Representative paths for the $\d_A$ and $\d_B$ generators.  $\d_A$ links $S^1_A$ (the blue circle) once in $A$ (the grey interior of the torus), and $\d_B$ links $S^1_B$ (the green circle) once in $B$ (the exterior of the torus).}
\label{fig5}
\end{figure}

It will be useful for our considerations in section \ref{sec:nested} to locate homotopic representative paths for $\d_A$ and $\d_B$.  For $\d_A$ we take such a path to be homotopic to one linking $S^1_A$ once, and similarly take $\d_B$ to link $S^1_B$ once.   It is easy to see that in terms of the coordinates introduced in section \ref{sec:definitions} and the choice of base points \eqref{specificP}, representatives of the $\d_A$ and $\d_B$ homotopy classes can be taken to be the circles
\begin{align}\label{deltaArep}
\d_A \simeq \{ \chi = \tfrac\pi4-\epsilon, \ \th=\tfrac\pi{nq}, \ 0\le\ps<2\pi \}.
\end{align}
Likewise, $\b$ is homotopic to $S^1_B$ which can be continuously deformed to the circle 
\begin{align}\label{deltaBrep}
\d_B \simeq \{ \chi = \tfrac\pi4+\epsilon, \ 0\le\th<2\pi,\ \ps=0 \}.
\end{align}
These representative paths for $\d_A$ and $\d_B$ are illustrated in figure \ref{fig5}.  We note, upon comparing to the description of the $\a$ and $\b$ paths given in \eqref{alpharep} and \eqref{betarep} (or in figure \ref{fig1}), that in $S^3\setminus S^1_A \cup S^1_B$, $\d_A$ is homotopic to $\b$ and $\d_B$ is homotopic to $\a$.

By choosing the deformation retracts in $\pi A$ and $\pi B$ to be the same as in the previous section --- i.e., by choosing only the products of $\ta_k$- and $\tb_k$-type elements in the groupoids $\pi A$ and $\pi B$ --- the group $F(x_0)$ defined by \eqref{fdef} and \eqref{Fdef} is the same as that constructed below \nameref{linkretract}.  In particular, it is independent of the new generators $\d_A$ and $\d_B$.

Since both $\pi C$ and $F(x_0)$ are unchanged from the last section, the relations with respect to which we take the quotient of the free product $\pi A(x_0) \ast \pi B(x_0) \ast F(x_0) $ are also unchanged.  Therefore, we have the same set of relations as in equation \eqref{conjugacyrelations}.
The fundamental group follows from \nameref{gsvkthm}, giving

\paragraph{Theorem 5.3}
\label{unknotpi1} 
\begin{align}\label{}
\pi X(x_0) &= \pi_1 (S^3 \setminus S^1_B \cup \kn \cup S^1_A) 
\nn\\
&= \Bigl\langle \a \, , \d_A \, ,  \b \, , \d_B \, , f_k \, \Big| \a \d_A \a^{-1}\d_A^{-1} \, , \b \d_B \b^{-1}\d_B^{-1} \, ,\, f_0 \, , 
\a^p f_k \b^{-q} f_k^{-1} \quad \forall k\in\Z_n \, \Bigr\rangle .
\nn
\end{align}

As a consistency check, \nameref{result1} of the previous section is  recovered by eliminating the additional generators by introducing the relations $\d_A = \I$ and $\d_B = \I$, which is equivalent to removing the two unknots from the initial setup.

\section{Knot group of nested torus links}
\label{sec:nested}

We now generalize to \emph{nested torus links}, $\fL_m$, $m\in\N$.  Conceptually, what we mean by this is a union of $m$ torus links, $\link{a}$, $a=1,\ldots,m$, each lying on a 2-torus, $T^{(a)}$, where the $T^{(a)}$ torus lies in the interior\footnote{This is the conventional --- not the topological --- use of ``interior'' as introduced in section \ref{sec:definitions}.} of $T^{(a+1)}$, and where $T^{(a+1)}$ is homotopic to $T^{(a)}$ within the interior of $T^{(a+1)}$.  


We can make this description more precise by defining this family of links inductively in the number of nestings, $m$.   When $m=1$, $\fL_1=\link{1}$ is a torus link.  In the coordinate system on $S^3$ introduced in section \ref{sec:definitions}, $\fL_1$ is on the torus $T$ at $\chi=\pi/4$, and so, in particular, all the links of $\fL_1$ are in $\bar{D_A}$, the closure of the solid donut, $D_A$, ``interior'' to $T$.  In coordinates, this is the set of all points in $S^3$ with $0\le\chi\le\pi/4$.  We will define nested torus links so that $\fL_m \subset \bar{D_A}$ for all $m$.

Consider the homeomorphism $\phi: S^3 \to S^3$ which acts on the $(\chi,\th,\ps)$ coordinates by $\phi: (\chi,\th,\ps) \mapsto (\frac2\pi \chi^2, \th, \ps)$.   This is easily checked to map the $0\le\chi\le\pi/2$ interval to itself and has the property on the open interval that $\phi(\chi) < \chi$.  Poetically, $\phi$ stretches and shrinks $S^3$ by compressing its foliation by constant-$\chi$ tori towards the $\chi=0$ ``interior'' circle $S^1_A$.  In particular, $\phi$ maps the $0\le\chi\le\pi/4$ interval to the $0\le\chi\le\pi/8$ interval, so $\phi(\bar{D_A})$ is strictly in the the topological interior of $\bar{D_A}$.

Given an $\fL_m$, define $\fL_{m+1}$ by
\begin{align}\label{fLdef}
\fL_{m+1} := \link{m+1} \cup \phi(\fL_m) ,
\end{align}
where, as defined in \eqref{Knqpdef}, $\link{m+1}$ lies on the $\chi=\pi/4$ torus.  It thus follows that $\fL_{m+1} \subset \bar{D_A}$.  It also follows that
\begin{align}\label{}
\fL_m &= \link{m} \cup \phi(\link{m-1}) \cup 
\phi^2(\link{m-2}) \cup \dots \cup \phi^{m-1}(\link{1}) .
\nn
\end{align}

We will compute the knot group of $\fL_m$ by computing the knot group of $S^1_B \cup \fL_m$ inductively in $m$ using the groupoid Seifert--van Kampen theorem.  The final result for $\fL_m$ then follows easily by ``filling in'' the $S^1_B$.

\paragraph{Lemma 6.1}
\label{lemma6.1} 

Let $X_1 := S^3 \minus S^1_B \cup \fL_1$ with $\fL_1\equiv \link{1}$.  Then
\begin{align}\label{}
\pi_1(X_1) 
&= \Bigl\langle \a_1 \, ,  \b_1 \, , \a_2 \, , f^{(1)}_{k_1} \, \Big|
\, \a_1^{p_1} f^{(1)}_{k_1} \b_1^{-q_1} f^{(1)-1}_{k_1} \, ,
\, \b_1 \a_2 \b_1^{-1} \a_2^{-1} \, , 
\, f^{(1)}_0 \, , 
\quad \forall k_1\in\Z_{n_1} 
\, \Bigr\rangle .
\nn
\end{align}
Furthermore, a representative path in the homotopy class $\a_2$ is that of $\d_B$ given in \eqref{deltaBrep}.

\noindent \textbf{Proof:}  
\nameref{unknotpi1} computes the knot group $\pi_1(X')$ with $X'= S^3 \minus S^1_B\cup\fL_1\cup S^1_A$.  To remove $S^1_A$ from this link, one ``fills in'' $S^1_A$ in $X'$.  That is, $X_1 = X' \cup S^1_A$.  Then, by an application of the Seifert--van Kampen theorem (we don't need the groupoid version in this case), we have that $\pi_1(X_1)$ is the amalgamated product (pushout of groups) $\pi_1(X_1) = \pi_1(X') \ast_{\pi_1(\til T)} \pi_1(S^1_A)$ where $\til T$ is homotopic to a torus which is the boundary of an arbitrarily small thickening of $S^1_A$.  Let $\ta$ be the generator of $\pi_1(S^1_A)\simeq\Z$.  Then the amalgamated product gives the presentation $\pi_1(X_1) = \left\langle \pi_1(X'), \ta \, \big| \, \ta \a^{-1}, \d_A \right\rangle$, where we are using the same notation for the generators of $\pi_1(X')$ as in \nameref{unknotpi1}.  But these relations amount to removing the additional $\ta$ generator, and adding a relation $\d_A=\I$ to those in \nameref{unknotpi1}.  Therefore $\pi_1(X_1)$ is as given in the statement of the lemma, up to a change in notation from \nameref{unknotpi1} where $\a$ and $\b$ are renamed $\a_1$ and $\b_1$, the $f_k$ are renamed $f^{(1)}_{k_1}$, and $\d_B$ is renamed $\a_2$ in preparation for the induction step.  
\hfill $\Box$

Note that the proof of this lemma justifies the claim made in the last sentence of section \ref{sec:withunknots} that imposing the relation $\d_A=\I$ is equivalent to removing the $S^1_A$ unknot from the link.

\paragraph{Theorem 6.2}
\label{nstdgrp}

Let $X_m := S^3 \minus S^1_B \cup \fL_m$. Then
\begin{align}
\pi_1 (X_m) 
& = \Big\langle 
\a_a ,  \b_a , \a_{m+1} , f^{(a)}_{k_a} \, \Big| \,
\a_a^{p_a} f^{(a)}_{k_a} \b_a^{-q_a} f^{(a) \, -1}_{k_a} ,\,  
f^{(a)}_0 ,\, 
\b_a \a_{a+1} \b_a^{-1}\a_{a+1}^{-1} ,
\nn\\
& \hspace{7cm}
\ \forall a \in\{1,\ldots,m\},\  \forall k_a\in\Z_{n_a}
\Big \rangle .\nn
\end{align}
Furthermore, $\d_B$ given in \eqref{deltaBrep} is a path in the homotopy class $\a_{m+1}$.

\noindent \textbf{Proof:}  By \nameref{lemma6.1} the theorem holds for $m=1$.  Assuming it is true for $m-1$, we will show it to hold for $m$, and then the theorem follows by induction.
We compute $\pi_1(X_m)$ given $\pi_1(X_{m-1})$ by applying \nameref{gsvkthm}.  Realize $X_m$ as the union $X_m=A\cup B$ where $A=\bar{D_A}\cap X_m$ and $B=\bar{D_B}\cap X_m$.  Recall that $\bar{D_A}$ is the closure of the ``interior'' of the $\chi=\pi/4$ torus $T$ while $\bar{D_B}$ is the closure of its ``exterior'', where the closure is defined relative to $S^3$.  By construction $X_m = S^3 \minus S^1_B \cup \link{m} \cup \phi(\fL_{m-1})$ where $S^1_B$ is exterior to $T$, $\link{m}$ is in $T$, and $\phi(\fL_{m-1})$ is interior to $T$.   Thus the topological interiors of $A$ and $B$ are $\text{int}(A) = D_A \minus \phi(\fL_{m-1})$, $\text{int}(B) = D_B \minus S^1_B$, and $C:= A\cap B = T \minus \link{m}$.  As in earlier sections, $A$ and $B$ are path connected while $C$ has $n_m$ path components.  Fix a base point set $P=\{x_0, \ldots, x_{n_m-1}\}$ as before with one point in each component of $C$.

To compute the pushout of the fundamental groupoids of the sets $A$, $B$, and $C$ defined above all with base point set $P$, we use their presentations
\begin{align}  
\pi A &= \Bigg\langle
\scalebox{0.75}{%
\begin{tikzpicture}[->,>=stealth',shorten >=1pt,thick,
    auto,node distance=2.0cm,
    main node/.style={circle, fill=black!20},
    baseline={([yshift=-.5ex]current bounding box.center)}]
\tikzset{every loop/.style={min distance=10mm,looseness=10}}
  \node[main node] (0) at (0,3) {$x_0$};
  \node[main node] (1) at (1.5,2.1) {$x_1$};
  \node (2) at (.5,.9) {};
  \node at (.35,.85) {$\bullet$};
  \node at (0,.8) {$\bullet$};
  \node at (-.35,.85) {$\bullet$};
  \node (n-2) at (-.5,.9) {};
  \node[main node] (n-1) at (-1.5,2.1) {$\phantom{x_0}$};
  \node at (-1.5,2.1) {\scriptsize $x_{n_m{-}1}$};
  \path
    (0) edge[in=45,out=135,loop,fill=red!60!black!30] node {} (0)
    (0) edge[bend left] node[above] {$\ta_1$} (1)
    (1) edge[bend left] node[above left] {$\ta_2$} (2)
    (n-2) edge[bend left] node[above right] {$\ta_{n_m{-}1}$} (n-1)
    (n-1) edge[bend left] node[above left] {$\ta_{n_m}$} (0);
  \node[red!60!black] at (2.3,4.3) {$\pi_1(X_{m-1})$};
  \draw[thick,red!60!black,->] (1.3,4.3) -- (0,4.0);
\end{tikzpicture} }
\hspace{-1cm} \Bigg| \ 
(\ta_{n_m}\cdots\ta_2\ta_1) \a_m^{-1}
\ \Bigg\rangle,
\nn\\
\pi B &= \Bigg\langle
\scalebox{0.75}{%
\begin{tikzpicture}[->,>=stealth',shorten >=1pt,thick,
    auto,node distance=2.0cm,
    main node/.style={circle,
    fill=black!20},
    baseline={([yshift=-.5ex]current bounding box.center)}]
\tikzset{every loop/.style={min distance=10mm,looseness=10}}
  \node[main node] (0) at (0,3) {$x_0$};
  \node[main node] (1) at (1.5,2.1) {$x_1$};
  \node (2) at (.5,.9) {};
  \node at (.35,.85) {$\bullet$};
  \node at (0,.8) {$\bullet$};
  \node at (-.35,.85) {$\bullet$};
  \node (n-2) at (-.5,.9) {};
  \node[main node] (n-1) at (-1.5,2.1) {$\phantom{x_0}$};
  \node at (-1.5,2.1) {\scriptsize $x_{n_m{-}1}$};
  \path
    (0) edge[loop above] node[above] {$\d_B$} (0)
    (0) edge[bend left] node[above] {$\tb_1$} (1)
    (1) edge[bend left] node[above left] {$\tb_2$} (2)
    (n-2) edge[bend left] node[above right] {$\tb_{n_m{-}1}$} (n-1)
    (n-1) edge[bend left] node[above] {$\tb_{n_m}$} (0);
\end{tikzpicture} }
\ \Bigg| \ 
(\tb_{n_m}\cdots\tb_2\tb_1) \d_B 
(\tb_{n_m}\cdots\tb_2\tb_1)^{-1} \d_B^{-1}
\ \Bigg\rangle ,
\nn\\
\pi C &= \left\langle
\scalebox{0.75}{%
\begin{tikzpicture}[->,>=stealth',shorten >=1pt,thick,
    auto,node distance=2.0cm,
    main node/.style={circle,
    fill=black!20},
    baseline={([yshift=-.5ex]current bounding box.center)}]
\tikzset{every loop/.style={min distance=10mm,looseness=10}}
  \node[main node] (0) at (0,3) {$x_0$};
  \node[main node] (1) at (1.5,2.1) {$x_1$};
  \node (2) at (.5,.9) {};
  \node at (.7,1.1) {$\bullet$};
  \node at (.25,.85) {$\bullet$};
  \node at (-.25,.85) {$\bullet$};
  \node at (-.7,1.1) {$\bullet$};
  \node (n-2) at (-.5,.9) {};
  \node[main node] (n-1) at (-1.5,2.1) {$\phantom{x_0}$};
  \node at (-1.5,2.1) {\scriptsize $x_{n_m{-}1}$};
  \path
    (0) edge[loop above,] node[above] {$\g_0$} (0)
    (1) edge[in=30,out=60,loop] node[above left] {$\g_1$} (1)
    (n-1) edge[in=120,out=150,loop] node[above] {$\g_{n_m-1}$} (n-1);
\end{tikzpicture} }
\right\rangle .\nn
\end{align} 
The presentations for $\pi B$ and $\pi C$ follow immediately from \nameref{lemma5.2} since $B$ and $C$ are the same as in that case.  In the presentation of $\pi A$ given above, the red-shaded loop based at $x_0$ represents the generators and relations of $\pi_1 (X_{m-1} := \pi_1(X_{m-1},x_0)$.  It appears because $A$ is homotopic to $X_{m-1}$.  This follows because $A$ is homotopic to $\text{int}(A)=D_A\minus \phi(\fL_{m-1})$, and $D_A$ is homeomorphic to $S^3\minus S^1_B$ under a mapping which takes $\chi\to2\chi$ with $\th$ and $\ps$ fixed.  Thus $A$ is homotopic to $S^3\minus S^1_B\cup\phi(\fL_{m-1}) = \phi(X_{m-1})$ since $\phi(S^1_B)=S^1_B$.  Because $A$ is connected, its object groups $\pi A(x_k)$ are all isomorphic, and so a generating graph of $\pi A$ consists of its object group $\pi A(x_0)$ together with a choice $\a_{x_k}$ for $k=1,\ldots,n_m-1$ of isomorphisms in $\pi A(x_0,x_k)$.  As in the definition of the retraction used in \nameref{linkretract}, choose the $\a_{x_k}:=\ta_k \dots \ta_2\ta_1$ for elements $\ta_k\in\pi A(x_{k-1},x_k)$ for  $k=1,\ldots,n_m-1$ defined similarly as in \nameref{linkgraph}.  Since $A$ is homotopic to $\phi(X_{m-1})$ which is homeomorphic to $X_{m-1}$, the object group $\pi A(x_0) = \pi_1(X_{m-1},x_0)$.  

Putting this all together gives the presentation for $\pi A$ shown above, but without the $\ta_{n_m}$ generator and without the additional relation shown.  By adding in $\ta_{n_m}$ as a generator, we generate an additional element, $\ta_{n_m}\dots\ta_2\ta_1$, in the object group at $x_0$.  But since $\pi_1(X_{m-1},x_0)$ is the full object group at $x_0$, this element must be one of the elements of $\pi_1(X_{m-1},x_0)$.  Indeed, this is the content of the additional relation shown which implies $\ta_{n_m}\dots\ta_2\ta_1 = \a_m$.  This particular relation follows from the induction hypothesis and the definition of the $\ta_k$.  In particular, the induction hypothesis states that $\a_m \in \pi_1(X_{m-1})$ is the homotopy class of the path $\d_B$ in \eqref{deltaBrep}.  Under the $\phi$ homeomorphism, this path is mapped to a similar circle but with $\chi \approx \frac\pi8+\epsilon$, so inside $D_A$.  On the other hand, from the definition of the $\ta_k$ given in \nameref{linkgraph}, it followed that the $\ta_{n_m}\dots\ta_2\ta_1$ loop is homotopic to the path $\a$ in \eqref{alpharep}.  Comparing the two, it is clear that they are homotopic by shifting the $\chi$ coordinate from $\frac\pi8+\epsilon$ to $\frac\pi4-\epsilon$.  This is a valid homotopy in $A$ since the region $\frac\pi8 < \chi <\frac\pi4$ is a subset of $A$: the excised $\phi(\fL_{m-1})$ link lies in the $\chi \le \frac\pi8$ region by construction.

With these presentations of the fundamental groupoids, their pushout is computed using \nameref{gsvkthm} in a similar way as it was in \nameref{result1} and \nameref{unknotpi1}.  In particular, choosing the deformation retraction $r$ along the $\a_{x_k}$ defined above and $s$ along $\b_{x_k}$ defined analogously, one finds that $ri\g_k = \a_m^{p_m}$ and $sj\g_k=\b_m^{q_m}$ for all $k\in\Z_{n_m}$ where $\b_m := \tb_{n_m} \dots \tb_2 \tb_1$.  Therefore the generators of $\pi_1(X_m)$ are those of $\pi_1(X_{m-1})$ together with $\b_m$, $f^{(m)}_{k_m}$ for $k_m\in\Z_{n_m}$, and $\d_B$.  And their relations are those of $\pi_1(X_{m-1})$, the usual ones involving the $f^{(m)}_{k_m}$ coming from \nameref{gsvkthm} --- namely, $f_0^{(m)}$ and $\a_m^{p_m} f_{k_m}^{(m)} \b_m^{-q_m} f_{k_m}^{(m)\,-1}$ --- and the relation $\b_m\d_B \b_m^{-1} \d_B^{-1}$ from $\pi B$.  Renaming $\d_B \to a_{m+1}$ then gives the presentation of $\pi_1(X_m)$ given in the statement of the theorem. 

\hfill $\Box$

An easy corollary is that the knot group of a nested torus link $\fL_m$ without the $S^1_B$ unknot is the same as that given in \nameref{nstdgrp} but with an additional relation setting $\a_{m+1}=\I$.  This follows from the same argument used in the proof of \nameref{lemma6.1}.  The previous results in this paper, \nameref{result1} and \nameref{unknotpi1}, all follow as special cases of \nameref{nstdgrp}.

\acknowledgments

It is a pleasure to thank Jacques Distler, Sergei Gukov, Cody Long, Matteo Lotito, Yongchao L\"u, Mario Martone, and Marko Stosic for helpful comments and discussions.  This work was supported in part by DOE grant DE-SC0011784.   DK was partially supported by a MUSE summer undergraduate research grant and Honors Student Award at University of Cincinnati.

\bibliographystyle{JHEP}

\providecommand{\href}[2]{#2}\begingroup\raggedright\endgroup

\end{document}